# Convex Regions and their 'Fairest' Equipartitioning Fans

R. Nandakumar (nandacumar@gmail.com)

**Abstract:** *A k-fan is a set of k half-lines (rays) all starting from the same point, called the origin of the fan. We discuss the partition of convex 2D regions into n (a positive integer) equal area convex pieces by fans with the following additional requirement: the perimeters of the resultant equal area pieces should be as close to one another as possible. We present some basic properties of such fans, which we call 'fairest equipartitioning fans', and raise further questions.*

## 1. Introduction: Problem Statement

A k-fan is a set of k half-lines (rays) all starting from the same point, which we call the origin of the fan. A fan is convex if successive rays are at an angle less than or equal to 180 degrees – this paper considers only convex fans. Given any convex region C in the plane and a number n (throughout our discussion C is a closed compact set formed by the union of its interior and its boundary). Let P be any point in the plane. If P is in the exterior of C, there is exactly one fan of n-1 rays with origin at P that partitions C into n equal area convex pieces. If P is in the interior of C, there could be infinitely many different n-fans originating at P that partition C into n equal area convex pieces.

**Note:** We consider a set of parallel lines as a fan originating at infinity.

**Definitions**

If point P is outside region C, there is exactly one fan with origin at P and breaking C into n equal area convex pieces. If P is in the interior, there could be infinitely many such fans yielding n equal area convex pieces. For each such convex equipartition of C into n pieces with fans, we define its **'fairness'** as the ratio between the maximum among perimeters of the pieces and the minimum among the piece perimeters. A fan partition for which fairness has the value 1, we call a **perfectly fair** fan (equi)partition (such a fan we call, for brevity, a 'perfect fan'). For each n, to each point on the plane P, we associate a unique value of a **'fairness function'** F(P,n) - if P is in the interior, F(P,n) is the fairness of that fan with **least** fairness ratio among all equipartitiong fans with origin at P. Now, for any C and n, with P ranging over the plane, we are interested in finding the global minimum of F(P,n) (its minimum among all fans originating at all points on the plane) – this minimum, we say, results in the **'fairest fan (equi)partition'** of C for that n. These definitions continue the definitions originally proposed in [1] and [2].

**Remarks**

1. It is obvious that for any n, a circle can be perfectly fair partitioned by infinitely many fans all originating at its center. Any parallelogram can also be perfectly fair partitioned for any n with equally spaced parallel rays. It **appears** that these two are the only convex shapes which allow perfect fan partitions for all n.

2. Any ellipse allows perfect fan partition for n=4 (fan origin at ellipse center) and for n=3 (there are two possible origins on the major axis on either side of the center). For no other n does any ellipse appear to allow perfectly fair fan equipartitions. It is also not clear which convex shape has the highest number of distinct perfectly fair fan origins for distinct values of n - just as the ellipse (probably) has 2.

For n=3, any C has a perfectly fair fan partition – this was proved in [3]. But, for a general n and general C, there may not be any perfectly fair fan partition. Hence the algorithmic problem of finding the fairest equipartitioning fan for a given C and n is of interest and most of our discussion will be on this problem.

3. The general 'fair partition conjecture' [1] states that for all n, any C allows partition (not necessarily partition by fans alone) into n pieces of equal area and equal perimeter. While open for general n, it has been established for all prime power values of n ([4], [5], [6]). Later, we show that partitioning by convex fans is much inferior to general fair partition in terms of maximizing fairness.

## 2. The n -> infinity limit of the Problem

The problem of finding the fairest equipartitioning fan (global minimum of the fairness function F(P,n)) for given C is easier to solve in the large n limit (n-> infinity). Indeed, for any point P on the plane (inside or outside C), we consider the infinite set of rays starting at P which have some intersection with C and its boundary at least at one point other than P (note: if P is in the exterior of C, a ray from P that is tangent to C at a vertex of C is **not** part of this set). Consider among such rays for a given P, the ratio between the maximum and minimum of their intersection lengths with C. It is easy to see that, in the large n limit, this ratio is also the ratio between the highest and lowest perimeter of equal area pieces separated by the fairest fan originating at P.

It can be seen that for any point P in the interior or on the boundary of C, F(P,n) as n->infinity has a finite value (note: starting in the deep interior, if P were to move close the boundary, the fairness function becomes very large and then discontinuously falls to a finite value as P hits the boundary). For most points outside C, F(P,n) tends to infinity with n; but it crucially tends to finite values for those $O(n^2)$ external points where edges of C extended outwards intersect (such intersection points could be infinitely far away, if two edges of C are parallel). At each such external intersection point, F(P,n), being finite there, obviously has a local minimum.

F(P,n) (n->infinity) has other local minima at (1) vertices of C, (2) mid points of edges of C and at a single point in the deep interior of C (obviously, this is the interior point where the ratio between the maximum distance from there to the boundary and the minimum distance to the boundary is least). All these minima together give a total of $O(n^2)$ candidate origins for the fairest equipartitioning fan of C for infinite n. Computing the fans and their fairness ratio for each external and boundary minimum is straightforward (for each point it is an O(m) calculation where m is the number of sides of C).

**Sub-problem:** Finding the single interior minimum of the fairness ratio (n-> infinity). As can be seen with simple examples, the origin of this fan need not lie on the medial axis of C. However, since there is only one such minimum in the interior for any C, it could be found by a suitable perturbative calculation.

Once all candidate fan origins and their fairness ratios are known, we select the minimum. Later in this paper (figure 2) we show an example illustrating the n->infinity limit.

## 3. Finite n

For finite n, finding the fairest fan of a given C is more difficult. We did following numerical analysis: for various test convex polygons, find F(P,n) at closely spaced points on the XY plane and plot these values along Z to see where the local and global minima lie on the resulting 'terrain' (when P is in the interior, we consider only a very large finite number of fans from the infinitely many possible ones so the calculations are close approximations).

The contour plots (we use gnuplot) for the fairness ratio (especially when C has symmetry) are often intriguingly beautiful. A set of specific examples are given below:

**Figure 1a.** C is a 12-sided polygon approximating an ellipse with semi-major axis 8 and semi-minor axis 5 and center at (0, 0) (the polygon is shown marked in black). We show F(P,n) for n=3. There are 2 minima, on the major axis, either side of the origin nearly 6 units away (not on the foci). The fairest fans originating there are indeed perfectly fair (F(P,n)= 1). For any ellipse and n=3, there are 2 such perfectly fair fan origins. (The analytical calculation of their positions has been posed as a challenge at [7]).

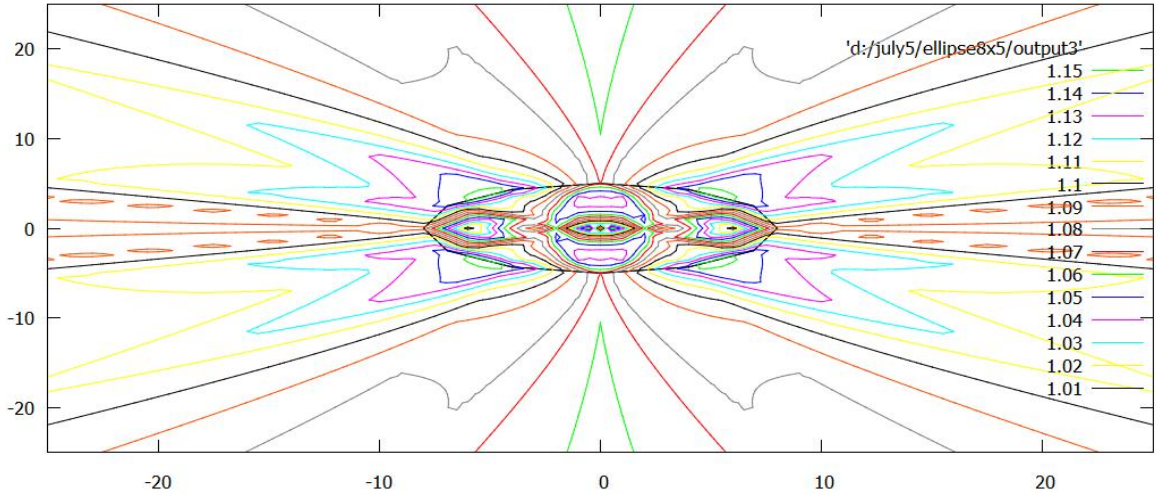

**Figure 1b:** the same polygon as figure 1a with n = 10. There are two minima of F(P,n) but they lie on the **minor** axis. The fairness ratio is slightly greater than 1, so the fairest fan partition is not perfectly fair.

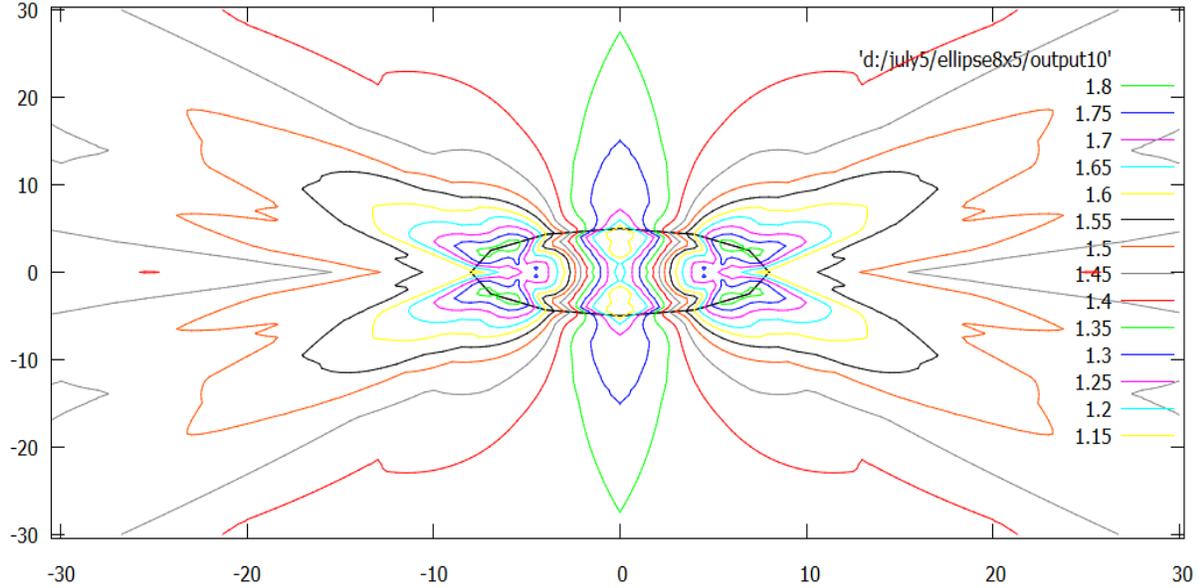

**Figure 1c:** the same polygon as above but with n = 100. There is now only a single minimum that lies at the center. The fairest partition is not perfectly fair.

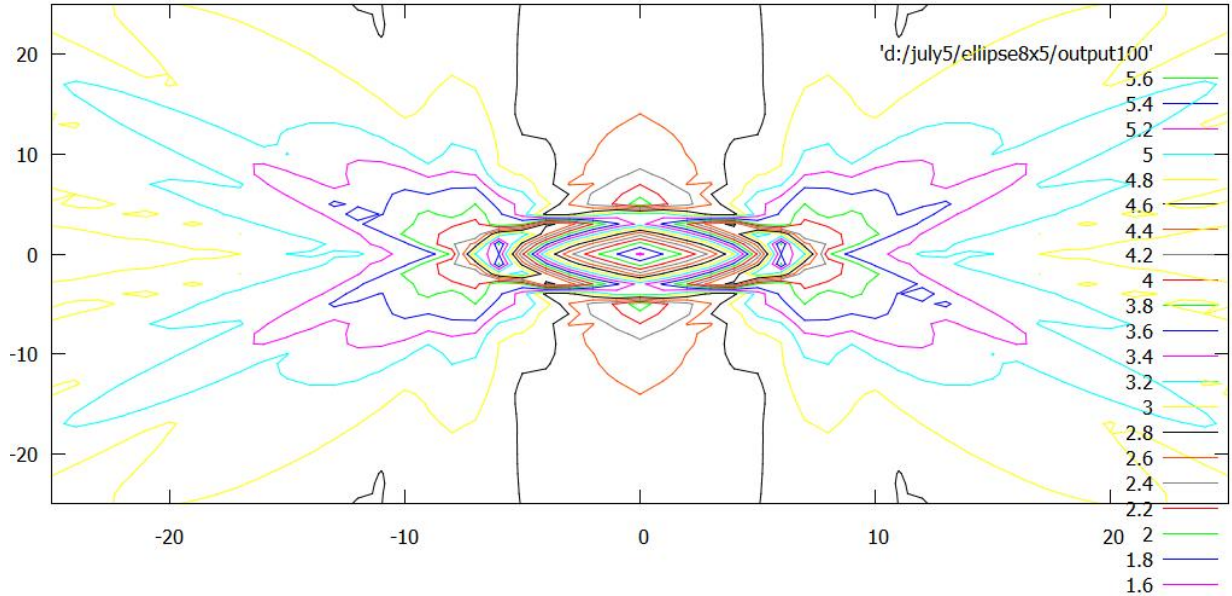

For further increment of n, the global minimum of F(P,n) stays at the origin. The 'terrain' keeps changing but the change with each increment of n is slight and the positions and alignment of the valleys and ridges stabilize. The minimum of fairness function stays greater than 1. As n tends to infinity, the fairest fan origin remains at the center and the fairness ratio there tends to the ratio between the major and minor axes of the ellipse.

For finite n and any general convex region C being fan partitioned, we found we find that there are two regimes – n very small and n large but finite. The following observations were made:

1. For small n, F(P,n) can have several local minima in the interior of C. The vertices of C need not be local minima or maxima. Then, for sufficiently large but finite n (typically of the order of 100), for any C, most of the interior is filled by a single large basin (which slowly grows with n and tends to an infinitesimally inset C in the large n limit) and this basin always has a single local minimum within it. In addition, for sufficiently large but finite n, the vertices and mid points of the edges develop local minima and basins around them; and these minima persist for larger values of n. (Note: Just inside C, all around its boundary runs a ridge of very high values of F).
2. In the exterior, near all the intersection points of extended edges the terrain develops basins (for sufficiently large but finite n) – for any n (finite or otherwise), each such basin always contains a single local minimum (not a set of closely spaced local minima).
3. As n is steadily incremented, the position of the global minimum of F(P,n) can jump around the plane (from one basin to another in the fairness function terrain; the basins themselves do not shift with n) and for sufficiently large but finite n, the minimum reaches the neighbourhood of the position for n=infinity and stays in that neighborhood for higher values of n (as an example, we have figures 1 a-c). As n tends to infinity, the external basins shrink in area to deep pits at the intersection points; everywhere else in the exterior, F(P,n) tends to infinity.
4. Crucially, for large but finite values of n, there are no other local minima anywhere on the plane other than those mentioned above. This means that for sufficiently large (but finite) n, the global minimum can be found by first calculating the $O(n^2)$ candidate points for n-> infinity and searching in the neighborhood of each by some gradient descent approximation method.
5. Far from the polygon, the terrain has no local minima for any n (small or large) and merely consists of ridges and valleys radiating away.
6. The local minima for finite n need not coincide with the n = infinity minima positions – they lie only within the neighbourhood of the n=infinity local minima. *Example:* For any triangle, the n=infinity local minima for the fairness function are obviously, its vertices and a single point in its interior. Later in this document, we show (figure 3) an isosceles triangle which, for n=6, has a local minimum of the fairness function near a vertex but just outside. For n=infinity, the local minimum in that neighborhood is of course, the vertex itself.

**Remark:** It is quite conceivable that for large but finite n, the above-mentioned method of approximation, descending from the n=infinity local minima should give the correct answer efficiently. For small values of n, with more than one local minimum in the interior, some more work might be needed.

For the large n behavior, we show in figure 2a and 2b below the n=700 fairness function terrain for the non-regular convex hexagon with vertices at (0,0),(10,0),(11,7),(1,12),(-4,10),(-4,4). Figure 2b is the close-up of the central portion of the region in figure 2a.

**Figure2a:**

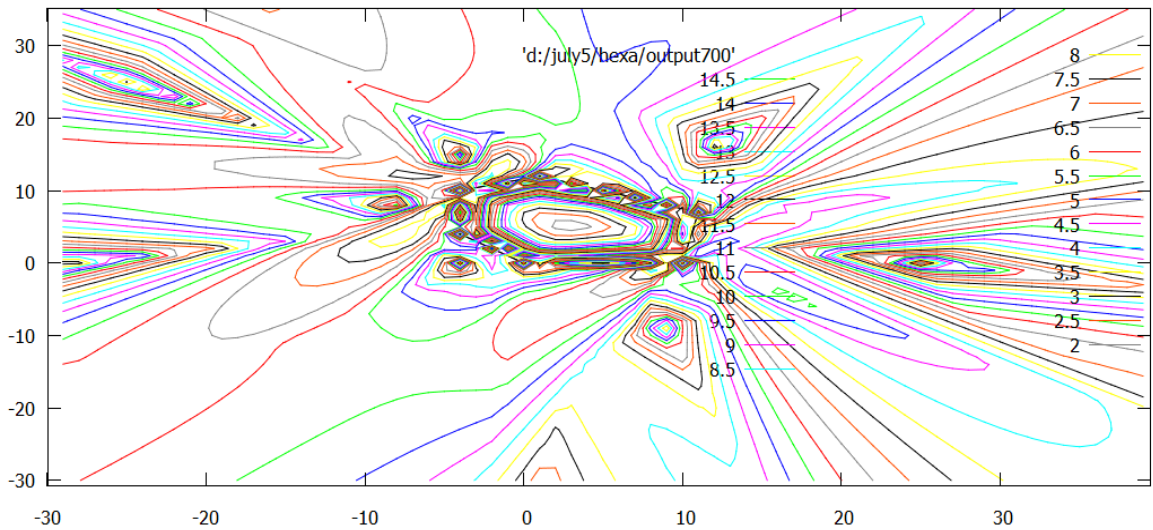

**Figure2b:**

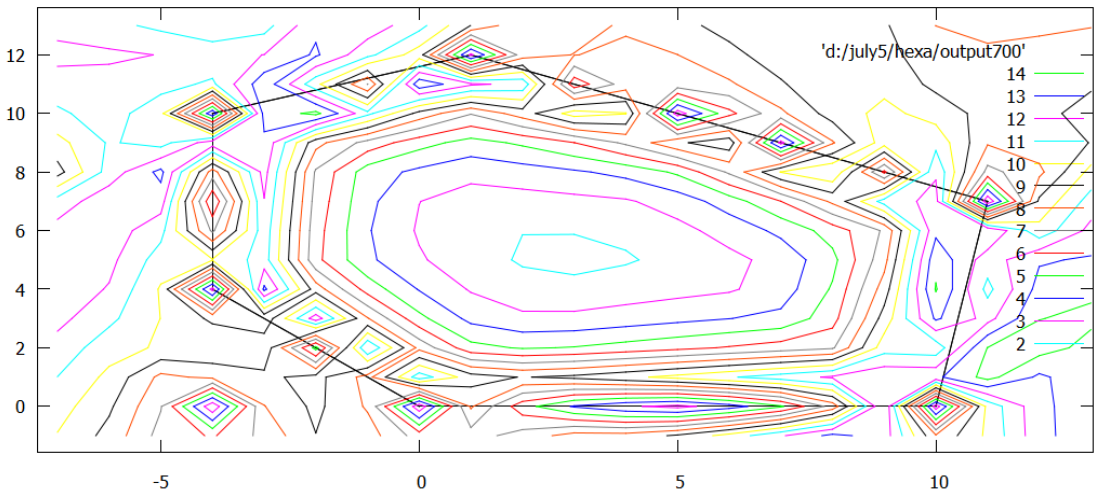

**Observations from Fig 2:** The terrain of the fairness function has a single basin in the interior of the polygon. There are local minima at the vertices and mid points of sides (the right top edge shows 2 local minima and not 1 but that is an artefact of the program, which samples the fairness function at a square grid of points). In the exterior, several basins are seen - they all have minima near the intersection points of the edges of the polygon extended outwards.

**Subproblem:** Consider finite n. As we know, with any single point P in the interior of C as origin, there may be infinitely many equipartitioning fans. How to analytically find the fan with the least fairness ratio?

*Answer sketch:* Consider any one convex fan equipartition of C (let C have m sides) into n pieces (numbered from 1 to n) with a convex fan originating at P; note the fairness ratio of this partition. Now, keeping fan origin at P, rotate this layout continuously (keeping the areas of pieces constant) until piece 1 coincides with the original position of piece 2, piece 2 with 3 and so forth. Obviously, during this rotation, the fairness ratio at that P would have changed continuously thru its entire possible range. Let us count the *changes in sign of the change in perimeter* of each piece during this rotation. It can be seen that the sum over all n pieces of this count of 'switches' would be only O(m). At any given stage of the rotation, the fairness ratio is determined only by the max and min among the piece perimeters. We could first do the perimeter and fairness ratio calculations for the n pieces only at the switches; then, changes in the max and min fairness ratios in the interval between any pair of successive switches would be straightforward to find because all piece perimeters change monotonically in such an interval. So, we infer that an algorithm taking polynomial time in m would suffice to find the least fairness ratio among fan partitions with a single origin.

**Further question:** Given a finite n: if we note the local minima for n= infinity (we know their positions) and evaluate the fairness function at only these points and choose the minimum (with no further perturbation calculations), how much will be off the best fairness ratio for that n? Can a bound (dependent on n) be found?

**Remarks:**

1. How does one upper-bound the number of local minima of F(P,n) for small n? Does it depend on the number of sides of C, on n itself and the symmetry of the polygon? The answer: The number of local minima does not appear to depend on n but only on the shape of C. The following behavior is seen for all test cases we studied: For low n (n within 10 or so), there usually are multiple local minima still O(n) in the interior. By intermediate values, there is only one local minimum inside the region but the vertices, midpoints of edges, external points of edge intersection. are not yet local minima - for sufficiently large n (n>100), for any polygon we checked, these special points would have become local minima and the interior filled by one big basin.

2. As n increases to large values, the fairest fan origin reaches and thereafter stays within the neighborhood of the fairest fan origin for n= infinity. But the finite value of n for which the best fan origin for that n reaches the neighborhood of the best fan center for n = infinity can be very large.

**Example:** Consider ellipses of high eccentricity. As noted earlier, its n->infinity fairest fan origin is the center of the ellipse. But our experiments clearly show that for up to some large value of n (which increases with eccentricity), the fairest fan origin lies on the boundary of the ellipse on its minor axis and only past a large critical n, it jumps to the center - to stay there thereafter on incrementing n.

3. We now compare fairest fan partitions and the fairest general convex partition of a convex region, the criterion for comparison being the difference between the least fairness ratios for the two partitions (as proved in [4], [5] and [6], this least fairness ratio is identically 1 for any C and any prime power value of n for *general* convex partitions).

**Proposition:** Equipartition with fans and the general fair partition problem are widely different in scope. Indeed, fan partitions cannot be good substitutes to general convex equipartitions as far as fairness (equality of perimeters of the resultant equal area pieces) goes - for there exist convex regions and values of n where the minimum of the fairness function of fairest convex partition by fans is arbitrarily greater (inferior) than for the fairest general partition.

**Supporting argument:** Above, we noted that for sufficiently large n, the fairest fan partition of **any** ellipse has the fan origin at its center. Consider n very large and also a prime power. For any such n, the fairness ratio for a fan originating at center of the ellipse is close to the ratio between the major and minor axes, which we can be arbitrarily large for an eccentric ellipse. Now, as proved in [4], [5] and [6], for n a prime power, the fairest general partition has fairness ratio exactly 1, so the above proposition is validated.

4. The above discussion on external local minima of the fairness function implies the following. If the boundary of convex region C has no straight line segment portions, there is no local minimum outside the polygon for n=infinity. We would like to ask further if for small and finite n, there could be external local minima for such shapes. To our knowledge, for such smooth convex regions, the local minima of the fairness function are all on or within the region. This also implies that for sufficiently large n, since the interior becomes a single basin for the F(P,n), such smooth regions have only one local minimum - which is close to the global minimum at n=infinity.

## 4. Further Questions and Conclusion

1. For a given n, what is that convex shape of specified area that is the **worst** region for fair partition by fans? By worst, we mean that minimum of the fairness function should be the **highest** possible. For large n, ellipses appear good candidates for the highest least fairness ratio for given area. At the other extreme, for n=4, the least fairness ratio for any ellipse is obviously a perfect 1 (fan origin at center) and our experiments show there do exist other convex shapes which have no perfectly fair fan partition for n=4 – we do not know which is the worst shape for n=4.

2. Which convex shape has largest number of distinct perfect fan origins for different n's? As we saw above, ellipses have distinct origins for n=3 (a pair of them) and n=4. For example, can one find a convex shape for which say there are m distinct points that work as origins for perfect fans for different n values? The circle apart, is there any shapes are there which allows infinitely many perfect fair partitions with fans originating not at infinity?

3. Considering only fans originating outside C but not at infinity, is there an upper bound on the n for which any C can be fan partitioned *perfectly*? Based on a suggestion from N. Ramana Rao, an isosceles triangle was found numerically that has a perfect fan with origin outside for n=6. Figure 3 shows a triangle with base vertices at (0, 0) and (10, 0) and apex at (5, 9.01042). Computations show a perfect fan origin for n =6 near the apex, a little outside the triangle at approximately (5, 9.6…).

**Figure3**

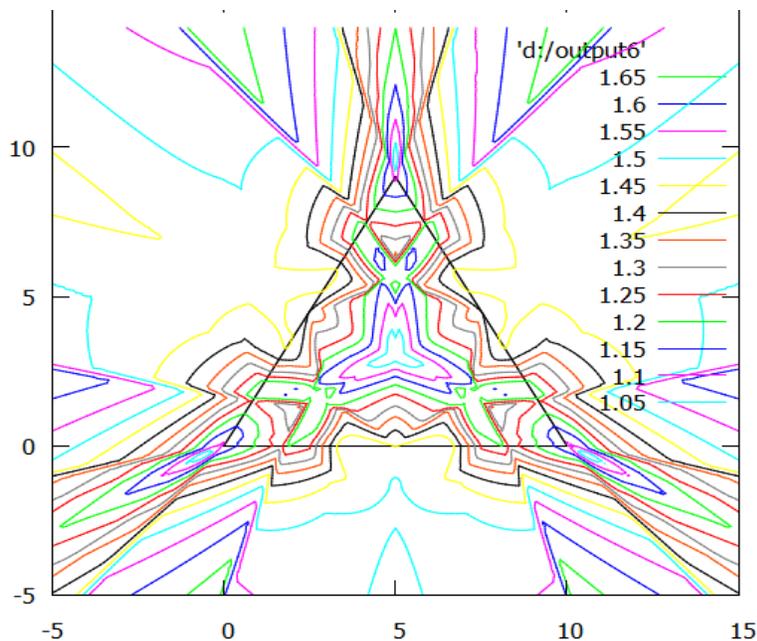

We could not find any isosceles triangle which has perfect fans with external origin for n = 8 and higher n. However, a convex region slightly deformed from an isosceles triangle might work at least for n=8. What we do not know now is "how precisely to deform?" It is quite likely that there is a finite n (and a corresponding convex region) beyond which no perfect partition from a fan with focus outside will work (the reason: there is no pair {a convex region, an outside point P} such that for n -> infinity, all lines from P and passing thru the region have same length segments inside the region). We do not know this critical n value- and what is the corresponding critical convex region.

4. Is there a convex region which has no perfectly fair fan partition for any n>3 (ellipses appear to have none for n>4)? A possibly related question: For any n, is there any convex region that can be fan partitioned with perfect fairness into n pieces that are all dissimilar in shape to one another? Note: For n=3, any convex equipartition is necessarily a fan partition and by [3], this question becomes trivial.

5. If we are to think of global maxima instead of minima of the fairness function: the question is trivial in the large n limit since most points on the plane have the fairness function tending to infinity. However for finite n and especially small n, the question could be of interest. More broadly, analytical properties of the fairness function could be explored further and that may help in finding the fairest fan for finite n.

6. Let a given convex region to be partitioned into n pieces with different areas: $a_1, a_2, \ldots a_n$. A possible additional specification: if $a_{min}$ and $a_{max}$ are to be the maximum and minimum of these piece areas, we could try to find fan partitions where ratio between $p_{min}$ and $p_{max}$ (the perimeters of the min and max area pieces respectively) is as close to the square root of $a_{min}/a_{max}$ as possible. Even in this case, our experimental setup woud generate unique contour patterns for each polygon and each such set of areas and many of these patterns could be visually interesting. This may be of some interest from the point of view of graphic design.

7. And finally, what about higher dimensions?

**Acknowledgements:** Thanks to N. Ramana Rao and John Rekesh for their generous help.